\documentclass{article}

\usepackage{amsmath, amsthm, amsfonts, amssymb}
\usepackage[all]{xy}

\theoremstyle{plain}
\newtheorem{thm}{Theorem}[section]

\newtheorem*{question}{Question}

\theoremstyle{definition}

\newcommand{\F}{\mathbb{F}}

\newcommand{\C}{\mathbb{C}}
\newcommand{\Aut}{\mathrm{Aut}}
\newcommand{\Inn}{\mathrm{Inn}}
\newcommand{\Out}{\mathrm{Out}}
\newcommand{\Hom}{\mathrm{Hom}}
\newcommand{\GL}{\mathrm{GL}}
\newcommand{\Sp}{\mathrm{Sp}}

\allowdisplaybreaks[1]

\title{A Survey on Automorphism Groups of Finite $p$-Groups}
\author{Geir T. Helleloid \\
				Department of Mathematics, Bldg. 380 \\
				Stanford University \\
				Stanford, CA 94305-2125 \\
				\texttt{geir@math.stanford.edu}
}

\begin{document}
\maketitle

\begin{abstract}
This survey on the automorphism groups of finite $p$-groups focuses on three major topics: explicit computations for familiar finite $p$-groups, such as the extraspecial $p$-groups and Sylow $p$-subgroups of Chevalley groups; constructing $p$-groups with specified automorphism groups; and the discovery of finite $p$-groups whose automorphism groups are or are not $p$-groups themselves.  The material is presented with varying levels of detail, with some of the examples given in complete detail.
\end{abstract}

\section{Introduction}
\label{introsec}

The goal of this survey is to communicate some of what is known about the automorphism groups of finite $p$-groups.  The focus is on three topics: explicit computations for familiar finite $p$-groups; constructing $p$-groups with specified automorphism groups; and the discovery of finite $p$-groups whose automorphism groups are or are not $p$-groups themselves.  Section~\ref{genthmsec} begins with some general theorems on automorphisms of finite $p$-groups.  Section~\ref{examplesec} continues with explicit examples of automorphism groups of finite $p$-groups found in the literature.  This includes the computations on the automorphism groups of the extraspecial $p$-groups (by Winter~\cite{win}), the Sylow $p$-subgroups of the Chevalley groups (by Gibbs~\cite{gib} and others), the Sylow $p$-subgroups of the symmetric group (by Bondarchuk~\cite{bon} and Lentoudis~\cite{len}), and some $p$-groups of maximal class and related $p$-groups.  Section~\ref{constructionsec} presents several theorems showing how to prescribe various quotients of the automorphism group of a finite $p$-group.  Section~\ref{ordersec} focuses on the order of the automorphism groups, concluding with many examples of finite $p$-groups whose automorphism group is a $p$-group.  Finally, Section~\ref{miscsec} contains some other miscellaneous results on the topic.  The material within is presented with varying levels of detail; in particular, some of the more explicit examples have been given in full detail.  Most of the necessary terminology is defined; some of the background material can be found by Huppert~\cite{hup}.

There are aspects of the research on the automorphism groups of finite $p$-groups that are largely omitted in this survey.  We mention three here.  The first is the conjecture that $|G| \le |\Aut(G)|$ for all non-cyclic finite $p$-groups $G$ of order at least $p^3$.  This has been verified for many families of $p$-groups, and no counter-examples are known; there is an old survey by Davitt~\cite{dav}.  The second is the (large) body of work on finer structural questions, like how the automorphism group of an abelian $p$-group splits or examples of finite $p$-groups whose automorphism group fixes all normal subgroups.  The third is the computational aspect of determining the automorphism group of a finite $p$-group.  Eick, Leedham-Green, and O'Brien~\cite{elo} describe an algorithm for constructing the automorphism group of a finite $p$-group.  This algorithm has been implemented by Eick and O'Brien in the GAP package AutPGroup~\cite{gap4}.  There are references in~\cite{elo} to other related research as well.

There are a few other survey papers that overlap with this one.  Corsi Tani~\cite{cor} has a survey of examples of finite $p$-groups whose automorphism group is a $p$-group; all these examples are included in the survey along with some others.  Starostin~\cite{sta} and Mann~\cite{man} have surveys of questions on finite $p$-groups, each of which includes a section on automorphism groups; Starostin focuses on specific examples related to the $|G| \le |\Aut(G)|$ conjecture and finer structural questions.

All $p$-groups mentioned in this survey will be finite and $p$ will always denote a prime.

\section{General Theorems}
\label{genthmsec}

Here we summarize some basic theorems about the automorphism group of a finite $p$-group $G$, for the most part following the survey of Mann~\cite{man}.  First, we can identify two subgroups of $\Aut(G)$ which are themselves $p$-groups.  Let $\Aut_c(G)$ be the automorphisms of $P$ which induce the identity automorphism on $G/Z(G)$ (these are called the \emph{central automorphisms} of $G$), and let $\Aut_f(G)$ be the automorphisms of $G$ which induces the identity automorphism on $G/\Phi(G)$.  Then $\Aut_c(G)$ and $\Aut_f(G)$ are $p$-groups.  More results on $\Aut_c(G)$ are given by Curran and McCaughan~\cite{cm}.

The next result is a theorem of Gasch\"{u}tz~\cite{gas}, which states that all finite $p$-groups have outer automorphisms.  Furthermore, unless $G \cong C_p$, there is an outer automorphism whose order is a power of $p$.  It is an open question of Berkovich as to whether this outer automorphism can be chosen to have order $p$.  Schmid~\cite{sch} extended Gasch\"{u}tz' theorem to show that if $G$ is a finite nonabelian $p$-group, then the outer automorphism can be chosen to act trivially on the center.  Furthermore, if $G$ is neither elementary abelian nor extra-special, then $\Out(G)$ has a non-trivial normal $p$-subgroup.  Webb~\cite{web3} proved Gash\"{u}tz's theorem and Schmid's first generalization in a simpler way and without group cohomology.  If $G$ is not elementary abelian or extra-special, then M\"{u}ller~\cite{mul} shows that $\Aut_f(G) > \Inn(G)$.

As mentioned in Section~\ref{introsec}, one prominent open question is whether or not  $|G| \le |\Aut(G)|$ for all non-cyclic $p$-groups $G$ of order at least $p^3$.  A related question concerns the automorphism tower of $G$, namely
\[
P_0 = G \to G_1 = \Aut(G) \to G_2 = \Aut(G_1) \to \cdots,
\]
where the maps are the natural maps from $G_i$ to $\Inn(G_i)$.  For general groups $G$, a theorem of Wielandt shows that if $G$ is centerless, then the automorphism tower of $G$ becomes stationary in a finite number of steps.  Little is known about the automorphism tower of finite $p$-groups.  In particular, it is not known whether there exist finite $p$-groups $G$ other than $D_8$ with $\Aut(G) \cong G$.

\section{The Automorphisms of Familiar $p$-Groups}
\label{examplesec}

There are several familiar families of finite $p$-groups whose automorphisms have been described in a reasonably complete manner.  The goal of this section to present these results as concretely as possible.  We begin with a nearly exact determination of the automorphism groups of the extraspecial $p$-groups.  The next subsection discusses the maximal unipotent subgroups of Chevalley groups and Steinberg groups, for which Gibbs~\cite{gib} describes six types of automorphisms that generate the automorphism group.  For type $A_l$, Pavlov~\cite{pav} and Weir~\cite{wei} have (essentially) computed the exact structure of the automorphism group.  The last three subsections summarize what is known about the automorphism groups of the Sylow $p$-subgroups of the symmetric group, $p$-groups of maximal class, and certain stem covers.  We note that Barghi and Ahmedy~\cite{ba} claim to determine the automorphism group of a class of special $p$-groups constructed by Verardi~\cite{ver}; unfortunately, as is pointed out in the MathSciNet review of~\cite{ba}, the proofs in this paper are incorrect.

\subsection{The Extraspecial $p$-Groups}

Winter~\cite{win} gives a nearly complete description of the automorphism group of an extraspecial $p$-group.  (Griess~\cite{gri} states many of these results without proof.)  Following Winter's exposition, we will present some basic facts about extraspecial $p$-groups and then describe their automorphisms.

Recall that a finite $p$-group $G$ is \emph{special} if either $G$ is elementary abelian or $Z(G) = G' = \Phi(G)$.  Furthermore, a non-abelian special $p$-group $G$ is \emph{extraspecial} if $Z(G) = G' = \Phi(G) \cong C_p$.  The order of an extraspecial $p$-group is always an odd power of $p$, and there are two isomorphism classes of extraspecial $p$-groups of order $p^{2n+1}$ for each prime $p$ and positive integer $n$, as proved in Gorenstein~\cite{gor}.  When $p = 2$, both isomorphism classes have exponent $4$.  When $p$ is odd, one of these isomorphism classes has exponent $p$ and the other has exponent $p^2$.

Any extraspecial $p$-group $G$ of order $p^{2n+1}$ has generators $x_1, x_2, \dots, x_{2n}$ satisfying the following relations, where $z$ is a fixed generator of $Z(G)$: 
\begin{alignat*}{3}
\left[x_{2i-1}, x_{2i}\right] &= z && \quad \textrm{for $1 \le i \le n$,} \\
\left[x_i, x_j\right] &= 1 && \quad \textrm{for $1 \le i,j \le n$ and $|i-j| > 1$, and} \\
x_i^p &\in Z(G) && \quad \textrm{for $1 \le i \le 2n$}.
\end{alignat*}

When $p$ is odd, either $x_i^p = 1$ for $1 \le i \le 2n$, in which case $G$ has exponent $p$, or $x_1^p = z$ and $x_i^p = 1$ for $2 \le i \le 2n$, in which case $G$ has exponent $p^2$.  When $p = 2$, either $x_i^2 = x_i^2 = 1$ for $1 \le i \le 2n$, or $x_1^2 = x_2^2 = z$ and $x_i^2 = 1$ for $3 \le i \le 2n$.

Recall that if two groups $A$ and $B$ have isomorphic centers $Z(A) \stackrel{\phi}{\cong} Z(B)$, then the \emph{central product} of $A$ and $B$ is the group
\[
(A \times B)/\{(z, z^{-1} \phi) \; : \; z \in Z(A)\}.
\]
All extraspecial $p$-groups can be written as iterated central products as follows.  If $p$ is odd, let $M$ be the extraspecial $p$-group of order $p^3$ and exponent $p$, and let $N$ be the extraspecial $p$-group of order $p^3$ and exponent $p^2$.  The extraspecial $p$-group of order $p^{2n+1}$ and exponent $p$ is the central product of $n$ copies of $M$, while the extraspecial $p$-group of order $p^{2n+1}$ and exponent $p^2$ is the central product of $n-1$ copies of $M$ and one copy of $N$.  If $p = 2$, the extraspecial $2$-group of order $2^{2n+1}$ and $x_1^2 = x_2^2 = 1$ is isomorphic to the central product of $n$ copies of the dihedral group $D_8$, while the extraspecial $2$-group of order $2^{2n+1}$ and $x_1^2 = x_2^2 = z$ is isomorphic to the central product of $n-1$ copies of $D_8$ and one copy of the quaternion group $Q_8$.

When $G$ has exponent $p$, we can view the group more concretely.  The extraspecial $p$-group of order $p^{2n+1}$ and exponent $p$ is isomorphic to the group of $(n+1) \times (n+1)$ matrices over $\F_p$ with ones along the diagonal, arbitrary entries in the rest of the first row and the last column, and zeroes elsewhere.

In~\cite{win}, Winter states the following theorem on the automorphism groups of the extraspecial $p$-groups for all primes $p$.

\begin{thm}[Winter~\cite{win}]
\label{spthm}
Let $G$ be an extraspecial $p$-group of order $p^{2n+1}$.  Let $I = \Inn(G)$ and let $H$ be the normal subgroup of $\Aut(G)$ which acts trivially on $Z(G)$.  Then
\begin{enumerate}
\item $I \cong (C_p)^{2n}$.
\item $\Aut(G) \cong H \rtimes \left< \theta \right>$, where $\theta$ has order $p-1$.
\item If $p$ is odd and $G$ has exponent $p$, then $H/I \cong \Sp(2n, \F_p)$, and the order of $H/I$ is $p^{n^2} \prod_{i=1}^n{(p^{2i}-1)}$.
\item If $p$ is odd and $G$ has exponent $p^2$, then $H/I \cong Q \rtimes \Sp(2n-2, \F_p)$, where $Q$ is a normal extraspecial $p$-group of order $p^{2n-1}$, and the order of $H/I$ is $p^{n^2} \prod_{i=1}^{n-1}{(p^{2i} - 1)}$.  The group $Q \rtimes \Sp(2n-2, \F_p)$ is isomorphic to the subgroup of $\Sp(2n, \F_p)$ consisting of elements whose matrix $(a_{ij})$ with respect to a fixed basis satisfies $a_{11} = 1$ and $a_{i1} = 0$ for $i > 1$.
\item If $p = 2$ and $G$ is isomorphic to the central product of $n$ copies of $D_8$, then $H/I$ is isomorphic to the orthogonal group of order $2^{n(n-1)+1}(2^n-1)\prod_{i=1}^{n-1}{(2^{2i}-1)}$ that preserves the quadratic form $\xi_1 \xi_2 + \xi_3 \xi_4 + \cdots + \xi_{2n-1} \xi_{2n}$ over $\F_2$.
\item If $p = 2$ and $G$ is isomorphic to the central product of $n-1$ copies of $D_8$ and one copy of $Q_8$, then $H/I$ is isomorphic to the orthogonal group of order $2^{n(n-1)+1}(2^n+1)\prod_{i=1}^{n-1}{(2^{2i}-1)}$ that preserves the quadratic form $\xi_1 \xi_2 + \xi_3 \xi_4 + \cdots + \xi_{2n-1}^2 + \xi_{2n-1} \xi_{2n} + \xi_{2n}^2$ over $\F_2$.
\end{enumerate}
\end{thm}

The automorphisms in $\Aut(G)$ can be described more explicitly.  First, to define the automorphism $\theta$ in Theorem~\ref{spthm}, let $m$ be a primitive root modulo $p$ with $0 < m < p$.  Define $\theta$ by $\theta(x_{2i-1}) = x_{2i-1}^m$ and $\theta(x_{2i}) = x_{2i}$ for $1 \le i \le n$ and by $\theta(z) = z^m$.

As for $\Inn(G)$, it is clear that $\Inn(G) \cong G/Z(G)$ acts trivially on $Z(G)$ and $G/Z(G)$ and that $|\Inn(G)| = p^{2n}$.  The elements of $\Inn(G)$ are given explicitly by the $p^{2n}$ automorphisms $\sigma$ where $\sigma(z) = z$ and  $\sigma(x_i) = x_i z^{d_i}$ for each $i$ and some integers $0 \le d_i < p$.

It remains to describe $H$.  For $x \in G$, let $\overline{x}$ denote the coset $Z(G) x$.  Now $G/Z(G)$ becomes a non-degenerate symplectic space over $\F_p$ with the symplectic form $(\overline{x}, \overline{y}) = a$, where $[x,y] = z^a$ and $0 \le a < p$.  The symplectic group $\Sp(2n, \F_p)$ acts on $G/Z(G)$, preserving the given symplectic form.  Let $T \in \Sp(2n, \F_p)$ and let $A = (a_{ij})$ be the matrix of $T$ relative to the basis $\{\overline{x}_i\}$ (with $0 \le a_{ij} < p$).  Each element $x \in G$ can be uniquely expressed as $x = \left( \prod_{i=1}^{2n}{x_i^{a_i}} \right) z^c$ with $0 \le a_i, c < p$.  Define $\phi : G \to G$ by
\[
\phi(x) = \left[ \prod_{i=1}^{2n} \left( \prod_{j=1}^{2n}{x_j^{a_{ij}}} \right)^{a_i} \right] z^c.
\]
Then $\phi$ is an automorphism of $G$ if and only if $T$ is in the subgroup of $\Sp(2n, \F_p)$ given in Theorem~\ref{spthm}.

While Winter's results do give a complete description of the automorphisms in $H$, we can say a bit more about the structure of $H$; namely, whether or not $H$ splits over $I$.  As Griess proves in~\cite{gri}, when $p = 2$, $H$ splits if $n \le 2$ and does not split if $n \ge 3$.  Griess also states, but does not prove, that when $p$ is odd, $H$ always splits over $I$.  This observation is also made in, and can be deduced from, Isaacs~\cite{isa2} and~\cite{isa3} and Glasby and Howlett~\cite{gh}.  A short exposition of this proof when $p$ is odd and $G$ has exponent $p$ was communicated via the group-pub-forum mailing list by Martin Isaacs~\cite{isa}.  Let $J/I$ be the central involution of the symplectic group $H/I$, and let $T$ be a Sylow 2-subgroup of $J$.  Then $J$ is normal in $H$, $|T| = 2$, and the non-identity element of $T$ acts on $I$ by sending each element to its inverse.  Then $1 = C_I(T) = I \cap N_H(I)$.  On the other hand, by the Frattini argument, $H = J N_H(T)$, and since $T \le J \cap N_H(T)$, it follows that $H = I N_H(T)$.  But this means that $N_H(T)$ is a complement of $I$ in $H$, and so $H$ splits over $I$.

According to Griess, the proof when $G$ has exponent $p^2$ is more technical.

\subsection{The Maximal Unipotent Subgroups of a Chevalley Group}

Associated to any simple Lie algebra $\mathcal{L}$ over $\C$ and any field $K$ is the \emph{Chevalley group $G$ of type $\mathcal{L}$ over $K$}.  Table~\ref{chevtab} lists the Chevalley groups of types $A_l$, $B_l$, $C_l$, and $D_l$ over the finite field $\F_q$, as given in Carter~\cite{car}.  A few clarifications are necessary: the entry for type $B_l$ requires that $\F_q$ have odd characteristic; $O_{2l+1}(\F_q)$ is the orthogonal group which leaves the quadratic form $\xi_1 \xi_2 + \xi_3 \xi_4 + \cdots + \xi_{2l-1} \xi_{2l} + \xi_{2l+1}^2$ invariant over $\F_q$; and $O_{2l}(\F_q)$ is the orthogonal group which leaves the quadratic form $\xi_1 \xi_2 + \xi_3 \xi_4 + \cdots + \xi_{2l-1} \xi_{2l}$ invariant over $\F_q$.  Gibbs~\cite{gib} examines the automorphisms of a maximal unipotent subgroup of a Chevalley group over a field of characteristic not two or three.  We are only interested in finite groups, so from now on we will let $K = \F_q$, where $\F_q$ has characteristic $p > 3$ and $q = p^n$.  After some preliminaries on maximal unipotent subgroups, we will present his results.

\begin{table}
\begin{center}
\begin{tabular}{|c|c|}
\hline
Type & Chevalley Group \\
\hline
$A_l$ & $\mathrm{PSL}_{l+1}(\F_q)$ \\
$B_l$ & $\mathrm{P(O_{2l+1}'(\F_q))}$ \\
$C_l$ & $\mathrm{PSp_{2l}(\F_q)}$ \\
$D_l$ & $\mathrm{P(O_{2l}'(\F_q))}$ \\
\hline
\end{tabular} \\
\parbox{3in}{\caption{\label{chevtab}The Chevalley groups of types $A_l$, $B_l$, $C_l$ and $D_l$.}}
\end{center}
\end{table}

Let $\Sigma$, $\Sigma^+$, and $\pi$ denote the sets of roots, positive roots, and fundamental roots, respectively, of $\mathcal{L}$ relative to some Cartan subalgebra.  Then the Chevalley group $G$ is generated by $\{x_r(t) \; : \; r \in \Sigma, t \in \F_q \}$.  One maximal unipotent subgroup $U$ of $G$ is constructed as follows.  As a set, 
\[
U = \{x_r(t) \; : \; r \in \Sigma^+, t \in \F_q \}.
\]
For any $r, s \in \Sigma^+$ and $t, u \in \F_q$, the multiplication in $U$ is given by
\begin{eqnarray*}
x_r(t) x_r(u) &=& x_r(t+u) \\
\left[x_s(u), x_r(t)\right] &=& \left\{
	\begin{array}{c@{\quad:\quad}l}
	1 & r+s \textrm{ is not a root} \\
	\prod\limits_{ir+js \in \Sigma}{x_{ir+js}(C_{ij,rs}(-t)^iu^j)} & r+s \textrm{ is a root}. \\
	\end{array} \right.
\end{eqnarray*}
Here $i$ and $j$ are positive integers and $C_{ij,rs}$ are certain integers which depend on $\mathcal{L}$.  The order of $U$ is $q^N$, where $N = |\Sigma^+|$, and $U$ is a Sylow $p$-subgroup of $G$.

Gibbs~\cite{gib} shows that $\Aut(G)$ is generated by six types of automorphism, namely graph automorphisms, diagonal automorphisms, field automorphisms, central automorphisms, extremal automorphisms, and inner automorphisms.  Let the subgroup of $\Aut(G)$ generated by each type of automorphism be denoted by $P$, $D$, $F$, $C$, $E$, and $I$ respectively.  Let $P_r$ be the additive group generated by the roots of $\mathcal{L}$ and let $r_N$ be the highest root.  Label the fundamental roots $r_1, r_2, \dots, r_l$.

\begin{enumerate}
\item \emph{Graph Automorphisms}:  An automorphism $\sigma$ of $P_r$ that permutes both  $\pi$ and $\Sigma$ induces a graph automorphism of $U$ by sending $x_r(t)$ to $x_r(t)$ for all $r \in \Sigma^+$ and $t \in \F_q$.  Graph automorphisms correspond to automorphisms of the Dynkin diagram, and so types $A_l \; (l > 1)$, $D_l \; (l > 4)$, and $E_6$ have a graph automorphism of order 2, while the graph automorphisms in type $D_4$ form a group isomorphic to $S_3$.
\item \emph{Diagonal Automorphisms}:  Every character $\chi$ of $P_r$ with values in $\F_q^{\ast}$ induces a diagonal automorphism which maps $x_r(t)$ to $x_r(\chi(r)t)$ for all $r \in \Sigma^+$ and $t \in \F_q$.
\item \emph{Field Automorphisms}:  Every automorphism $\tau$ of $\F_q$ induces a field automorphism of $U$ which maps $x_r(t)$ to $x_r(\sigma(t))$ for all $r \in \Sigma^+$ and $t \in \F_q$.
\item \emph{Central Automorphisms}:  Let $\tau_i$ be endomorphisms of $\F_q^+$.  These induce a central automorphism that maps $x_{r_i}(t)$ to $x_{r_i}(t) x_{r_N}(\sigma_i(t))$ for $i = 1, \dots, l$ and all $t \in \F_q$.
\item \emph{Extremal Automorphisms}:  Suppose $r_j$ is a fundamental root such that $r_N - r_i$ is also a root.  Let $u \in \F_q^{\ast}$.  This determines an extremal automorphism which acts trivially on $x_{r_i}(t)$ for $i \neq j$ and sends $x_{r_j}(t)$ to 
\[
x_{r_j}(t) x_{r_N-r_j}(ut) x_{r_N}((1/2)N_{r_N-r_j,r_j}ut^2).
\]
Here, $N_{r_N-r_j, r_j}$ is a certain constant that depends on the type.  In type $C_l$, $r_N - 2r_j$ is also a root, and the map that acts trivially on $x_i(t)$ for $i \neq j$ and sends $x_{r_j}(t)$ to
\[
x_{r_j}(t) x_{r_N-2r_j}(ut) x_{r_N-r_j}((1/2)N_{r_N-2r_j,r_j}ut^2) x_{r_N}((1/3)C_{12,r_N-r_j,r_j}ut^3)
\]
is also an automorphism of $U$.
\end{enumerate}

Steinberg~\cite{ste} showed that the automorphism group of a Chevalley group of a finite field is generated by graph, diagonal, field, and inner automorphisms, which shows that $P$, $D$, and $F$ are, in fact, subgroups of $\Aut(U)$.  It is easy to see that the central automorphisms are automorphisms, and a quick computation verifies this for the extremal automorphisms as well.  Note that multiplying an extremal automorphism by a judicious choice of central automorphism, the $x_{r_N}(\cdot)$ term in the description of the extremal automorphisms disappears.  Therefore, a functionally equivalent definition of extremal automorphisms omits the $x_{r_N}(\cdot)$ term, and this is what we will use for what follows.

Gibbs does not compute the precise structure of $\Aut(U)$.  This has been done in type $A_l$ however; Pavlov~\cite{pav} computes $\Aut(U)$ over $\F_p$, while Weir~\cite{wei} computes it over $\F_q$ (although his computations contain a mistake which we will address in a moment).  We will present the result for type $A_l$ as explicitly as possible, pausing to note that it does seem feasible to compute the structure of $\Aut(U)$ for other types in a similar manner.

As mentioned before, in type $A_l$, we can view $U$ as the set of $(l+1) \times (l+1)$ upper triangular matrices with ones on the diagonal and arbitrary entries from $\F_q$ above the diagonal.  There are $\binom{l+1}{2}$ positive roots in type $A_l$, given by $r_i + r_{i+1} + \cdots + r_j$ for $1 \le i \le j \le l$.  Let $E_{i,j}$ be the $(l+1) \times (l+1)$ matrix with a $1$ in the $(i,j)$-entry and zeroes elsewhere.  Then $x_r(t) = I + tE_{i,j}$, where $r = r_i + r_{i+1} + \cdots + r_j$.  In particular, $x_{r_i}(t) = I + tE_{i,i+1}$.  In type $A_l$, some of the given types of automorphisms admit simpler descriptions; we will be content to describe their action on the elements $x_{r_i}(t)$.

As mentioned, in type $A_l$ there is one nontrivial graph automorphism of order $2$, and it acts by reflecting all matrices in $U$ across the anti-diagonal.  The diagonal automorphisms correspond to selecting $\chi_1, \dots, \chi_n \in \F_q^{\ast}$ and mapping $x_{r_i}(t)$ to $x_{r_i}(\chi_i t)$.  This is equivalent to conjugation by a diagonal matrix of determinant 1.  The diagonal automorphisms form an elementary abelian subgroup of order $(p-1)^{n-1}$.  The field automorphisms of $\F_q$ are generated by the Frobenius automorphism and form a cyclic subgroup of order $n$.

Let $a_1, \dots, a_n$ generate the additive group of $\F_q$.  Then $x_{r_i}(a_j)$ generate $U$.  The central automorphisms are generated by the automorphisms $\tau_j^m$ which send $x_{r_j}(a_m)$ to $x_{r_j}(a_m) x_{r_N}(1)$ and fix $x_{r_i}(a_k)$ for $i \neq j$ and $k \neq m$, where $j = 2, \dots, l-1$ and $m = 1, \dots, n$.  (When $j = 1$ or $j = l$, this automorphism is inner.)  The extremal automorphisms are generated by the automorphism that sends $x_{r_1}(t)$ to $x_{r_1}(t) x_{r_N-r_1}(t)$ and fixes $x_{r_i}(t)$ for $2 \le i \le l$ and the automorphism that sends $x_{r_l}(t)$ to $x_{r_l}(t) x_{r_N-r_l}(t)$ and fixes $x_{r_i}(t)$ for $1 \le i \le l-1$.  Finally the inner automorphism group is, of course, isomorphic to $U/Z(U)$ (the center of $U$ is generated by $x_{r_N}(t)$).

It is not hard to use these descriptions to deduce that 
\[
\Aut(U) \cong ((I \rtimes (E \times C)) \rtimes (D \rtimes F)) \rtimes P.
\]
Furthermore $E \times C$ is elementary abelian of order $q^{n(l-2)+2}$, $D$ is elementary abelian of order $(q-1)^l$, $F$ is cyclic of order $n$, and $I$ has order $q^{(l^2+l-2)/2}$.  It follows that the order of $\Aut(U)$ is
\[
2 n (q-1)^l q^{(l^2+l+2nl-4n+2)/2}.
\]

The error in Weir's paper~\cite{wei} stems from his claim that any $g \in \GL_n(\F_p)$ acting on $\F_q$ induces an automorphism of $U$ that maps $x_{r_i}(a_k)$ to $x_{r_i}(g(a_k))$, generalizing the field automorphisms.  However, it is clear that $g$ must, in fact, be a field automorphism, as for any $t, u \in \F_q$, 
\[
\left[ x_{r_1}(t), x_{r_{2}}(u) \right] = x_{r_1+r_2}(tu) = \left[ x_{r_1}(tu), x_{r_2}(1) \right].
\]
Applying $g$ to all terms shows that $g(tu) = g(t)g(u)$.

\subsection{Sylow $p$-Subgroups of the Symmetric Group}

The automorphism groups of Sylow $p$-subgroups of the symmetric group for $p > 2$ were examined independently by Bondarchuk~\cite{bon} and Lentoudis~\cite{len, len2, len3, lt}.  Their results are reasonably technical.  They do show that the order of the automorphism group of the Sylow $p$-subgroup of $S_{p^m}$, which is isomorphic to the $m$-fold iterated wreath product of $C_p$, is
\[
(p-1)^m p^{n(m)},
\]
where
\[
n(m) = p^{m-1} + p^{m-2} + \cdots + p^2 + \frac{1}{2} (m^2-m+2)p - 1.
\]

\subsection{$p$-Groups of Maximal Class}

A $p$-group of order $p^n$ is of \emph{maximal class} if it has nilpotence class $n - 1$.  It is not too hard to prove some basic results about the automorphism group of an arbitrary $p$-group of maximal class.  Our presentation follows Baartmans and Woeppel~\cite[Section 1]{bw}.

\begin{thm}
Let $G$ be a $p$-group of maximal class of order $p^n$, where $n \ge 4$ and $p$ is odd.  Then $\Aut(G)$ has a normal Sylow $p$-subgroup $P$ and $P$ has a $p'$-complement $H$, so that $\Aut(G) \cong H \rtimes P$.  Furthermore, $H$ is isomorphic to a subgroup of $C_{p-1} \times C_{p-1}$.
\end{thm}

The proof of this theorem begins by observing that $G$ has a characteristic cyclic series $G = G_0 \lhd G_1 \lhd \cdots \lhd G_n = 1$; that is, each $G_i$ is characteristic and $G_i/G_{i+1}$ is cyclic (see Huppert~\cite[Lemmas 14.2 and 14.4]{hup}).  By a result of Durbin and McDonald~\cite{dm}, $\Aut(G)$ is supersolvable and so has a normal Sylow $p$-subgroup $P$ with $p'$-complement $H$, and the exponent of $\Aut(G)$ divides $p^t (p-1)$ for some $t > 0$.  The additional result about the structure of $H$ comes from examining the actions of $H$ on the characteristic cyclic series and on $G/\Phi(G)$.  Baartmans and Woeppel remark that the above theorem holds for any finite $p$-group $G$ with a characteristic cyclic series.

Baartmans and Woeppel~\cite{bw} follow up these general results by focusing on automorphisms of $p$-groups of maximal class of exponent $p$ with a maximal subgroup which is abelian.  More specifically, the characteristic cyclic series can be taken to be a composition series, in which case $G_i = \gamma_i(G)$ for $i \ge 2$ and $G_1 = C_G(G_2/G_4)$, where $\gamma_i(G)$ is the $i$-th term in the lower central series of $G$.   Baartmans and Woeppel assume that $G_2$ is abelian.

In this case, they show by construction that $H \cong C_{p-1} \times C_{p-1}$.  Furthermore, $P$ is metabelian of nilpotence class $n-2$ and order $p^{2n-3}$.  (Recall that a \emph{metabelian group} is a group whose commutator subgroup is abelian.)  The subgroup $\Inn(G) \cong G/Z(G)$ has order $p^{n-1}$ and maximal class $n-2$.  The commutator subgroup $P'$ is the subgroup of $\Inn(G)$ induced by $G_2$.  Baartmans and Woeppel do explicitly describe the automorphisms of $G$, but the descriptions are too complicated to include here.

Other authors who investigate automorphisms of certain finite $p$-groups of maximal class include: Abbasi~\cite{abb}; Miech~\cite{mie}, who focuses on metabelian groups of maximal class; and Wolf~\cite{wol}, who looks at the centralizer of $Z(G)$ in certain subgroups of $\Aut(G)$.

Finally, in~\cite{juh}, Juh\'{a}sz considers more general $p$-groups than $p$-groups of maximal class.  Specifically, he looks at $p$-groups $G$ of nilpotence class $n-1$ in which $\gamma_1(G)/\gamma_2(G) \cong C_{p^m} \times C_{p^m}$ and $\gamma_i(G)/\gamma_{i+1}(G) \cong C_{p^m}$ for $2 \le i \le n-1$.  He refers to such groups as being of type $(n,m)$.  Groups of type $(n,1)$ are the $p$-groups of maximal class of order $p^n$.

Assume that $n \ge 4$ and $p > 2$.  As with groups of maximal class, the automorphism group of a group $G$ of type $(n,m)$ is a semi-direct product of a normal Sylow $p$-subgroup $P$ and its $p'$-complement $H$, and $H$ is isomorphic to a subgroup of $C_{p-1} \times C_{p-1}$.  Juh\'{a}sz' results are largely technical, dealing with the structure of $P$, especially when $G$ is metabelian.

\subsection{Stem Covers of an Elementary Abelian $p$-Group}
\label{stemcoversec}

In~\cite{web2}, Webb looks at the automorphism groups of stem covers of elementary abelian $p$-groups.  We start with some preliminaries on stem covers.  A group $G$ is a \emph{central extension} of $Q$ by $N$ if $N$ is a normal subgroup of $G$ lying in $Z(G)$ and $G/N \cong Q$.  If $N$ lies in $[G, G]$ as well, then $G$ is a \emph{stem extension} of $Q$.  The Schur multiplier $M(Q)$ of $Q$ is defined as the second cohomology group $H^2(Q, \C^{\ast})$, and it turns out that $N$ is isomorphic to a subgroup of $M(Q)$.  Alternatively, $M(Q)$ can be defined as the maximum group $N$ so that there exists a stem extension of $Q$ by $N$.  Such a stem extension is called a $\emph{stem cover}$.

Webb takes $Q$ to be elementary abelian of order $p^n$ with $p$ odd and $n \ge 2$.  Let $G$ be a stem cover of $Q$.  Then $N = Z(G) = [G,G] = M(Q) \cong Q \wedge Q$ and has order $p^{\binom{n}{2}}$.  Therefore $\Aut_c(G)$ are the automorphisms of $G$ which act trivially on $G/N \cong Q$.  Each automorphism $\alpha \in \Aut_c(G)$ corresponds uniquely to a homomorphism $\overline{\alpha} \in \Hom(Q, N)$ via the relationship $g N \overline{\alpha} = g^{-1} \cdot g\alpha$ for all $g \in G$.  Of course, $\Hom(Q, N)$ is an elementary abelian $p$-group of order $n \binom{n}{2}$, and so $\Aut(G)$ is an extension of a subgroup of $\Aut(Q) \cong \GL(n, \F_p)$ by an elementary abelian $p$-group of order $n \binom{n}{2}$.  Webb proves that the subgroup of $\Aut(Q)$ in question is usually trivial, leading to her main theorem.

\begin{thm}[Webb~\cite{web2}]
Let $G$ be elementary abelian of order $p^n$ with $p$ odd.  As $n \to \infty$, the proportion of stem covers of $G$ with elementary abelian automorphism group of order $p^{n \binom{n}{2}}$ tends to 1.
\end{thm}

\section{Quotients of Automorphism Groups}
\label{constructionsec}

Not every finite $p$-group is the automorphism group of a finite $p$-group.  A recent paper in this vein is by Cutolo, Smith, and Wiegold~\cite{csw}, who show that the only $p$-group of maximal class which is the automorphism group of a finite $p$-group is $D_8$.  But there are several extant results which show that certain quotients of the automorphism group can be arbitrary.

\subsection{The Central Quotient of the Automorphism Group}

\begin{thm}[Heineken and Liebeck~\cite{hl2}]
\label{hlthm}
Let $K$ be a finite group and let $p$ be an odd prime.  There exists a finite $p$-group $G$ of class $2$ and exponent $p^2$ such that $\Aut(G)/\Aut_c(G) \cong K$.
\end{thm}

The construction given by Heineken and Liebeck can be described rather easily.  Let $K$ be a group on $d$ generators $x_1, x_2, \dots, x_d$.  Let $D'(K)$ be the directed Cayley graph of $K$ relative to the given generators.  Form a new digraph $D(K)$ by replacing every arc in $D'(K)$ by a directed path of length $i$ if the original arc corresponded to the generator $x_i$.  Then $\Aut(D(K)) = K$.

Let $v_1, v_2, \dots, v_m$ be the vertices of $D(K)$.  Let $G$ be the $p$-group generated by elements $v_1, v_2, \dots, v_m$ where
\begin{enumerate}
\item $G'$ is the elementary abelian $p$-group freely generated by
\[
\{ [v_i, v_j] \; : \; 1 \le i < j \le m \}.
\]
\item For each vertex $v_i$, if $v_i$ has outgoing arcs to $v_{i_1}, v_{i_2}, \dots, v_{i_k}$, then
\[
v_i^p = [v_i, v_{i_1} \cdots v_{i_k}].
\]
\end{enumerate}
Heineken and Liebeck show that $\Aut(G)/\Aut_c(G) \cong K$ when $|K| \ge 5$.  (They give a special construction for $|K| < 5$.)  As Webb~\cite{web2} notes, $G$ is a special $p$-group.

They are actually able to determine the automorphism group of $G$ much more precisely, at least when $|K| \ge 5$.  Let the vertices of $D'(K)$ be called \emph{group-points}; they are naturally identified with vertices of $D(K)$.  Let $S$ be the set of vertices of $D(K)$ consisting of the group-point $e$ corresponding to the identity of $K$ and all vertices that can be reached along a directed path from $e$ that does not pass through any other group-points.  Assume that the vertices of $D(K)$ are labeled so that $v_1, \dots, v_s$ are the elements of $S$.  The central automorphisms which fix $v_{s+1}, v_{s+2}, \dots, v_m$ generate an elementary abelian $p$-group $U$ of rank $s |G'| = (1/2)ks^2(ks-1)$.  Every central automorphism of $G$ is of the form $\prod_{v \in K}{v^{-1} \alpha_v v}$, where the elements $\alpha_v \in U$ and $v \in K$ are uniquely determined.  Thus $\Aut_c(G)$ is the direct product of the conjugates of $U$ in $\Aut(G)$ and $\Aut(G) = U \wr K$.  It follows that $\Aut(G)$ has order $k p^l$, where $l = (1/2)k^2 s^2(ks-1)$ and $s = (1/2)d(d+1) + 1$ when $d \ge 2$ and $s = 1$ when $d = 1$.

Lawton~\cite{law} modified Heineken and Liebeck's techniques to construct smaller groups $G$ with $\Aut(G)/\Aut_c(G) \cong K$.  He uses undirected graphs which are much smaller, and the $p$-group $G$ which he defines is significantly simpler.

Webb~\cite{web} uses similar, though more complicated techniques, to obtain further results.  She defines a class of graphs called \emph{$Z$-graphs}; it turns out that almost all finite graphs are $Z$-graphs (that is, the proportion of graphs on $n$ vertices which are $Z$-graphs goes to 1 as $n$ goes to infinity).  To each $Z$-graph $\Lambda$, Webb associates a special $p$-group $G$ for which $\Aut(G)/\Aut_c(G) \cong \Aut(\Lambda)$.  The set of all special $p$-groups that arise from $Z$-graphs on $n$ vertices is denoted by $\mathcal{G}(p, n)$.  

\begin{thm}[Webb]
Let $p$ be any prime.  Then almost all of the groups in $\mathcal{G}(p,n)$ (as $n \to \infty$) have automorphism group $(C_p)^r$, where $r = n^2(n-1)/2$.
\end{thm}

The reason the group $(C_p)^r$ arises as the automorphism group is that for $G \in \mathcal{G}(p, n)$, $\Aut_c(G)$ is isomorphic to $\Hom(G/Z(G), Z(G))$, and hence to $(C_p)^r$ (for essentially the same reason as in Subsection~\ref{stemcoversec}).  Webb then shows that $\Aut(G)/\Aut_c(G)$ is usually trivial.

\begin{thm}[Webb]
\label{webbthm}
Let $K$ be a finite group which is not cyclic of order at most five.  Then for any prime $p$, there is a special $p$-group $G \in \mathcal{G}(p, 2|K|)$ with $\Aut(G)/\Aut_c(G) \cong K$.
\end{thm}

In particular, Theorem~\ref{webbthm} extends Heineken and Liebeck's result to the case $p = 2$.  Note that in Theorems~\ref{hlthm} and~\ref{webbthm}, the constructed groups are special and $\Aut_c(G) = \Aut_f(G)$, so that these theorems also prescribe $\Aut(G)/\Aut_f(G)$.

The $p = 2$ analogue of Heineken and Liebeck's result was discussed by Hughes~\cite{hug}.

\subsection{The Quotient $\Aut(G)/\Aut_f(G)$}

Bryant and Kov\'{a}cs~\cite{bk} look at prescribing the quotient $\Aut(G)/\Aut_f(G)$, taking a different approach from that of Heineken and Liebeck in that they assign $\Aut(G)/\Aut_f(G)$ as a linear group (and they do not bound the class of $G$).

\begin{thm}[Bryant and Kov\'{a}cs~\cite{bk}]
Let $p$ be any prime.  Let $K$ be a finite group with dimension $d \ge 2$ as a linear group over $\F_p$.  Then there exists a finite ($d$-generator) $p$-group $G$ such that $\Aut(G)/\Aut_f(G) \cong K$.
\end{thm}

This theorem is (essentially) non-constructive, in contrast to the results of Heineken and Liebeck.  To understand the main idea, let $F$ be the free group on $d$ generators and let $F_n$ be the $n$-th term in the Frattini series of $F$.  There is an action of $\GL(d, \F_p)$ on $F_n/F_{n+1}$.  If $U$ is a normal subgroup of $F$ with $F_{n+1} \le U \le F_n$, then $G = F/U$ is a finite $p$-group and $\Aut(G)/\Aut_f(G)$ is isomorphic to the normalizer of $U$ in $\GL(d,\F_p)$ (see~\cite[Theorems 2.7 and 2.8]{hm} for more details).  Bryant and Kov\'{a}cs show that if $n$ is large enough, then $F_n/F_{n+1}$ contains a regular $\F_p \GL(d, \F_p)$-module, which shows that any subgroup $K$ of $\GL(d,\F_p)$ occurs as the normalizer of some normal subgroup $U$ of $F$ with $F_n \le U \le F_{n+1}$.

\section{Orders of Automorphism Groups}
\label{ordersec}

The first two subsections describe some general theorems about the orders of automorphism groups of finite $p$-groups.  The third subsection gives the order of the automorphism group of an abelian $p$-group, and the last subsection offers many explicit examples of $p$-groups whose automorphism group is a $p$-group.  Helleloid and Martin~\cite{hm} have proved that, in several asymptotic senses, the automorphism group of a finite $p$-group is almost always a $p$-group.  However, a question raised in Mann~\cite{man} remains unanswered:

\begin{question}
Fix a prime $p$.  Let $v_n$ be the proportion of $p$-groups with order at most
$p^n$ whose automorphism group is a $p$-group.  Is it true that $\lim_{n \to \infty}{v_n} = 1$?
\end{question}

\subsection{Nilpotent Automorphism Groups}

In~\cite{yin}, Ying states two results about the occurrence of automorphism groups of $p$-groups which are $p$-groups, the second being a generalization of a result of Heineken and Liebeck~\cite{hl}.

\begin{thm}
If $G$ is a finite $p$-group and $\Aut(G)$ is nilpotent, then either $G$ is cyclic or $\Aut(G)$ is a $p$-group.
\end{thm}

\begin{thm}
Let $p$ be an odd prime and let $G$ be a finite two-generator $p$-group with cyclic commutator subgroup.  Then $\Aut(G)$ is not a $p$-group if and only if $G$ is the semi-direct product of an abelian subgroup by a cyclic subgroup.
\end{thm}

Heineken and Liebeck~\cite{hl} also have a criterion which determines whether or not a two-generator $p$-group of class 2 has an automorphism of order 2 or if the automorphism group is a $p$-group.  If $p$ is an odd prime and $G$ is a $p$-group that admits an automorphism which inverts some non-trivial element of $G$, then $G$ is an \emph{s.i. group} (a some-inversion group).  Clearly if $G$ is an s.i. group, it has an automorphism of order 2.  If $G$ is not an s.i. group, it is called an \emph{n.i. group} (a no-inversion group).

\begin{thm}
Let $p$ be an odd prime and let $G$ be a two-generator $p$-group of class 2.  Choose generators $x$ and $y$ such that
\[
\left< x, G' \right> \cap \left< y, G' \right> = G',
\]
and suppose that
\[
\left< x \right> \cap G' = \left< x^{p^m} \right>, \qquad \left< y \right> \cap G' = \left< y^{p^n} \right>.
\]
\begin{enumerate}
\item If either $x^{p^m} = 1$ or $y^{p^n} = 1$, then $G$ is an s.i. group.
\item If $x^{p^m} = [x,y]^{rp^k} \neq 1$ and $y^{p^n} = [x,y]^{sp^l} \neq 1$ with $(r,p) = (s,p) = 1$, and $(n-l+k-m)(k-l)$ is non-negative, then $G$ is an s.i. group.
\item If $m$, $n$, $k$, and $l$ are defined as in (2) and $(n-l+k-m)(k-l)$ is negative, then $G$ is an n.i. group and its automorphism group is a $p$-group
\end{enumerate}
\end{thm}

\subsection{Wreath Products}

In~\cite{hor2}, Horo{\v{s}}evski{\u\i} gives the following two theorems on the order of the automorphism group of a wreath product.

\begin{thm}
Let $A$ and $B$ be non-identity finite groups, and let $A_1$ be a maximal abelian subgroup of $A$ which can be distinguished as a direct factor of $A$.  Then
\[
\pi(\Aut(A \wr B)) = \pi(A) \cup \pi(B) \cup \pi(\Aut(A)) \cup \pi(\Aut(B)) \cup \pi(\Aut(A_1 \wr B)).
\]
\end{thm}

\begin{thm}
Let $P_1, P_2, \dots, P_m$ be non-identity finite $p$-groups.  Then
\[
\pi(\Aut(P_1 \wr P_2 \wr \cdots \wr P_m)) = \bigcup_{i=1}^m{\pi(\Aut(P_i))} \cup \{p\}.
\]
\end{thm}

Thus given any finite $p$-groups whose automorphism group is a $p$-group, we can construct infinitely many more by taking iterated wreath products.

\subsection{The Automorphism Group of an Abelian $p$-Group}

Macdonald~\cite[Chapter II, Theorem 1.6]{mac} calculates the order of the automorphism group of an abelian $p$-group.  The literature does contain some more technical results on the structure of such an automorphism group.

\begin{thm}
Let $G$ be an abelian $p$-group of type $\lambda$.  Then
\[
|\Aut(G)| = p^{|\lambda| + 2n(\lambda)} \prod_{i \ge 1}{\phi_{m_i(\lambda)}(p^{-1})},
\]
where $m_i(\lambda)$ is the number of parts of $\lambda$ equal to $i$, $n(\lambda) = \sum_{i \ge 1}{\binom{\lambda_i'}{2}}$, and $\phi_m(t) = (1-t)(1-t^2) \cdots (1-t^m)$.
\end{thm}

\subsection{Miscellaneous $p$-Groups Whose Automorphism Group is a $p$-Group}

In this subsection, we collect constructions of finite $p$-groups whose automorphism groups are $p$-groups.

The first example of a finite $p$-group whose automorphism group is a $p$-group was given by Miller~\cite{mil}, who constructed a non-abelian group of order 64 with an abelian automorphism group of order 128.  Generalized Miller's construction, Struik~\cite{str} gave the following infinite family of 2-groups whose automorphism groups are abelian $2$-groups:
\begin{eqnarray*}
G &=& \left< a,b,c,d \; : \; a^{2^n} = b^2 = c^2 = d^2 = 1, \right. \\
&& \qquad \left. [a,c] = [a,d] = [b,c] = [c,d] = 1, bab = a^{2^{n-1}}, bdb = cd \right>,
\end{eqnarray*}
where $n \ge 3$.  ($G$ can be expressed as a semi-direct product as well.)  Struik shows that $\Aut(G) \cong (C_2)^6 \times C_{2^{n-2}}$.  (As noted in~\cite{str}, it turns out that revision problem \#46 on p. 237 of Macdonald~\cite{mac2} asks the reader to show that $\Aut(G)$ is an abelian $2$-group.)  Also, Jamali~\cite{jam} has constructed, for $m \ge 2$ and $n \ge 3$, a non-abelian $n$-generator group of order $2^{2n+m-2}$ with exponent $2^m$ and abelian automorphism group $(C_2)^{n^2} \times C_{2^{m-2}}$.

More examples of 2-groups whose automorphism groups are 2-groups are given by Newman and O'Brien~\cite{no}. As an outgrowth of their computations on $2$-groups of order dividing 128, they present (without proof) three infinite families of 2-groups for which $|G| = |\Aut(G)|$.  They are, for $n \ge 3$,
\begin{enumerate}
\item $C_{2^{n-1}} \times C_2$,
\item $\left< a, b \; : a^{2^{n-1}} = b^2 = 1, a^b = a^{1+2^{n-2}} \right>$, and
\item $\left< a, b, c \; : \; a^{2^{n-2}} = b^2 = c^2 = [b, a] = 1, a^c = a^{1+2^{n-4}}, b^c = ba^{2^{n-3}} \right>$.
\end{enumerate}

Moving on to finite $p$-groups where $p$ is odd, for each $n \ge 2$ Horo{\v{s}}evski{\u\i}~\cite{hor} constructs a $p$-group with nilpotence class $n$ whose automorphism group is a $p$-group, and for each $d \ge 3$ he constructs a $p$-group on $d$ generators for each $d \ge 3$ whose automorphism group is a $p$-group.  (He gives explicit presentations for these groups.)

Curran~\cite{cur2} shows that if $(p-1,3) = 1$, then there is exactly one group of order $p^5$ whose automorphism group is a $p$-group (and it has order $p^6$).  It has the following presentation:
\begin{eqnarray*}
G &=& \left< a,b \; : \; b^p = [a,b]^p = [a,b,b]^p = [a,b,b,b]^p = [a,b,b,b,b] = 1, \right. \\
&& \qquad \left. a^p = [a,b,b,b] = [b,a,b]^{-1} \right>.
\end{eqnarray*}
When $(p-1,3) = 3$, there are no groups of order $p^5$ whose automorphism group is a $p$-group.  However, in this case, there are three groups of order $p^5$ which have no automorphisms of order 2.  Curran also shows that $p^6$ is the smallest order of a $p$-group which can occur as an automorphism group (when $p$ is odd).

Then, in~\cite{cur3}, Curran constructs 3-groups $G$ of order $3^n$ with $n \ge 6$ where $|\Aut(G)| = 3^{n+3}$ and $p$-groups $G$ for certain primes $p > 3$ with $|\Aut(G)| = p|G|$.  The MathSciNet review of~\cite{cur3} remarks that F. Menegazzo notes that for odd $p$ and $n \ge 3$, the automorphism group of 
\[
G = \left< a,b \; : \; a^{p^n} = 1, b^{p^n} = a^{p^{n-1}}, a^b = a^{1+p} \right>
\]
has order $p |G|$.

Ban and Yu~\cite{by2} prove the existence of a group $G$ of order $p^n$ with $|\Aut(G)| = p^{n+1}$, for $p > 2$ and $n \ge 6$.  In~\cite{hl}, Heineken and Liebeck construct a $p$-group of order $p^6$ and exponent $p^2$ for each odd prime $p$ which has an automorphism group of order $p^{10}$.

Jonah and Konvisser~\cite{jk} exhibit $p+1$ nonisomorphic groups of order $p^8$ with elementary abelian automorphism group of order $p^{16}$ for each prime $p$. All of these groups have elementary abelian and isomorphic commutator subgroups and commutator quotient groups, and they are nilpotent of class two. All their automorphisms are central.

Malone~\cite{mal} gives more examples of $p$-groups in which all automorphisms are central: for each odd prime $p$, he constructs a nonabelian finite $p$-group $G$ with a nonabelian automorphism group which comprises only central automorphisms.  Moreoever, his proof shows that if $F$ is any nonabelian finite $p$-group with $F' =Z(F)$ and $\Aut_c(F)=\Aut(F)$, then the direct product of $F$ with a cyclic group of order $p$ has the required property for $G$.

Caranti and Scoppola~\cite{cs} show that for every prime $p > 3$, if $n \ge 6$, there is a metabelian $p$-group of maximal class of order $p^n$ which has automorphism group of order $p^{2(n-2)}$, and if $n \ge 7$, there is a metabelian $p$-group of maximal class of order $p^n$ with an automorphism group of order $p^{2(n-2)+1}$.  They also show the existence of non-metabelian $p$-groups ($p > 3$) of maximal class whose automorphism groups have orders $p^7$ and $p^9$.

\section{Miscellaneous Results}
\label{miscsec}

This section contains a brief mention of several results which seem to be worth including.

In~\cite{ay}, Adney and Yen examine the automorphism group of a finite $p$-group $G$ of class 2 with no abelian direct factor, where $p$ is odd.  Under certain conditions on $G$, they show that $|G|$ divides $|\Aut(G)|$; the MathSciNet review of this article states that M. Newman can prove this with no extra conditions on $G$.  Ban and Yu have several papers on which groups can be the automorphism group of a $p$-group, focusing on groups of small order.  As an example, see~\cite{by}.  

Beisiegel~\cite{bei} shows that if $G$ is a $p$-group and not elementary abelian, then $\Aut(G)$ is $p$-constrained.  Furthermore, if $G$ has a cyclic commutator subgroup and $p > 2$, then $\Aut(G)$ is not a $p$-group if and only if $\Aut(G)$ contains an involution.  Menegazzo~\cite{men} studies the automorphism groups of finite non-abelian 2-generated $p$-groups with cyclic commutator subgroup for odd primes $p$.  He exhibits presentations for the relevant groups and computes the orders of $\Aut(G)$, $O_p(\Aut(G))$, and $\Aut_f(G)$.

Liebeck~\cite{lie} obtains an upper bound for the class of $\Aut_f(G)$.  Finally, Wang and Zhang~\cite{wz} discuss the automorphism groups of some $p$-groups.

\bibliographystyle{amsplain}
\bibliography{grouptheory}

\end{document}